\documentclass[12pt]{article}
\usepackage{amssymb,amsfonts,amsmath,geometry,esint}

\usepackage{color}

\def\R{\mathbb{R}}

\def\f{\varphi}

\def\irn{\int\limits_{\R^n}}

\def\Dshalf{\left(-\Delta\right)^{\!\frac {s}{2}}\!}  

\def\sstar{{2^*_s}}

\def\QED{\hfill {$\square$}\goodbreak \medskip}

\newtheorem{Theorem}{Theorem}
{Lemma}
\newtheorem{Proposition}
{Proposition}
{Corollary}
\newtheorem{Remark}
{Remark}

\date{}

\begin{document}

\title 
{A note on higher order \\fractional Hardy-Sobolev inequalities}

\author{Roberta Musina\footnote{Dipartimento di Scienze Matematiche, Informatiche e Fisiche, Universit\`a di Udine,
via delle Scienze, 206 -- 33100 Udine, Italy. Email: {roberta.musina@uniud.it}. 
Partially supported by PRID project {\em VARPROGE}.
}~ 
and \setcounter{footnote}{6}
Alexander I. Nazarov\footnote{
St.Petersburg Dept of Steklov Institute, Fontanka 27, St.Petersburg, 191023, Russia, 
and St.Petersburg State University, 
Universitetskii pr. 28, St.Petersburg, 198504, Russia. E-mail: al.il.nazarov@gmail.com.
Partially 
supported by RFBR grant 20-01-00630.
}
}

\date{}

\maketitle

\begin{abstract}
 We establish some qualitative properties of minimizers in the fractional Hardy--Sobolev inequalities of arbitrary order.
\end{abstract}


Assume $n\ge 1$ is a given integer, and take exponents  $s\in(0,\frac{n}{2})$, $q$, $ \beta$ satisfying
\begin{equation*}
2<q<\sstar:=\frac {2n}{n-2s}~,\qquad 
\frac {n}{q}-\beta=\frac {n}{2}-s~\!.
\end{equation*}
By H\"older interpolation between the  
Hardy and Sobolev inequalities one obtains the existence of a positive best constant $S_q$ such that
\begin{equation}
\label{eq:HS_ineq}
S_q\cdot\||x|^{-\beta} u\|_q^2\le\|\Dshalf u\|_2^2,\qquad u\in\mathcal D^s(\R^n),
\end{equation}
where 
$$
\mathcal D^s(\R^n)=\big\{u\in L^\sstar(\R^n)~|~\Dshalf u \in L^2(\R^n)~\!\big\}~,\quad 
{\mathcal F}\big[\Dshalf u\big] = |\xi|^{s}{\mathcal F}[u]
$$
and $\mathcal F$ is the Fourier transform in $\R^n$.

The existence of a minimizer  for $S_q$ was proved in \cite{Yang}. 
Also, it was claimed in \cite[Theorem 1.2]{Yang}, that any minimizer of (\ref{eq:HS_ineq}) has constant sign, and
is radially 
symmetric and strictly decreasing. However, the proof of this statement contains serious gaps. 
Namely, author of \cite{Yang} refers to the results of \cite{SV} and \cite{FS}, which were  formulated and proved only for $s<1$ 
and which are in general false for $s>1$, compare with \cite{MNTr}. 
Moreover, even in the case $s<1$ the argument in \cite{Yang} does not provide the {\it strict} monotonicity of the 
minimizer\footnote{In fact, for $s<1$ this can be proved by adapting the moving plane argument in \cite{CLO} 
or \cite{DMPS}, as  pointed out in \cite{MN20}.}.

In this note we prove the following statement.

\begin{Theorem}
\label{T:1}
Let $u$ be a minimizer for $S_q$. Up to a change of sign, $u$ is everywhere positive, radially symmetric and
strictly monotone decreasing with respect to $|x|$.
\end{Theorem}

\medskip
Before proving Theorem \ref{T:1} we recall some  notation. For $f\ge 0$ measurable and vanishing at infinity, $f^*$ stands for the symmetric-decreasing rearrangement of $f$, see \cite[Section 3.2]{LiLo}.

Accordingly with \cite{Cho}, \cite{FV}, for nonnegative measurable functions  $f,g$ vanishing at infinity, we write $f\prec g$ if 
$$
\int\limits_{B_r(0)}f^*(x)\,dx\le\int\limits_{B_r(0)}g^*(x)\,dx\quad \text{for any $r>0$}.
$$

Trivially, the pointwise inequality $f\le g$ implies $f\prec g$, while the inverse implication is not true, in general. Finally, we recall the next result.

\begin{Proposition} {\rm (see, e.g., \cite[Lemma 2.1]{FV})}.
The relation $f\prec g$ is equivalent to any of the following statements:
\begin{enumerate}
 \item 
 For any non-negative convex functions $\phi$ such that $\phi(0)=0$, we have
\begin{equation}
\label{eq:prec1}
\irn \phi(f^*(x))\,dx\le\irn \phi(g^*(x))\,dx~\!.
\end{equation}


\item
For any non-negative symmetric-decreasing function $\f$ we have
\begin{equation}
\label{eq:prec}
\irn f^*(x)\varphi(x)\,dx\le\irn g^*(x)\varphi(x)\,dx~\!;
\end{equation}
\end{enumerate}
\end{Proposition}

\paragraph{Proof of Theorem \ref{T:1}.}
We introduce the functions $U, V\in \mathcal D^s(\R^n)$ as (unique) weak solutions to 
\begin{equation}
\label{eq:s/2}
\Dshalf U= f:=\big|\Dshalf u\big|\in L^2(\R^n),\qquad
\Dshalf V= f^*\in L^2(\R^n),
\end{equation}
respectively. It is well known that $U$ and $V$ are explicitly given by convolutions
\begin{equation}
\label{eq:conv}
U(x)=C(n,s)\irn \frac {f(y)}{|x-y|^{n-s}}\,dy, \qquad V(x)=C(n,s)\irn \frac {f^*(y)}{|x-y|^{n-s}}\,dy~\!,
\end{equation}
where $C(n,s)$ is an explicitly known constant. In particular, $U, V>0$.

The Riesz rearrangement inequality, see, e.g., \cite[Sec. 3.7]{LiLo}, immediately implies (\ref{eq:prec}), 
and thus $U\prec V$, or, equivalently, $U^*\prec V^*=V$.

Further, we proceed in 3 steps.
\medskip

\noindent
{\bf Step 1}. We claim that $u$ strictly preserves the sign.
Indeed,  (\ref{eq:s/2}) and (\ref{eq:conv}) give
$$
 \Dshalf (U\pm u)=\big|\Dshalf u\big|\pm \Dshalf u\ge 0\quad
 \Rightarrow\quad U\pm u\ge0 \quad \Rightarrow\quad U\ge |u|.
 $$
Therefore,
 $$
 \||x|^{-\beta} U\|_q^2\ge \||x|^{-\beta} u\|_q^2 \quad\mbox{and}\quad \|\Dshalf U\|_2^2=\|\Dshalf u\|_2^2.
$$
Thus $U$ achieves $S_q$ as well. But then 
$$
\||x|^{-\beta} U\|_q^2= \||x|^{-\beta} u\|_q^2\quad \Longrightarrow\quad 
|u|\equiv U>0,
$$
and the claim follows. From now on, we assume that $u>0$.
\medskip

\noindent
{\bf Step 2}. We claim that $u=u^*$. Indeed, by Step 1 we have $u^*=U^*\prec V$. If we assume that $u\neq u^*$, then a basic rearrangement inequality, see, e.g., \cite[Sec. 3.4]{LiLo} together with  (\ref{eq:prec1}), (\ref{eq:prec}) gives
$$
\irn|x|^{-\beta q}u^q~\!dx < \irn|x|^{-\beta q}(u^*)^q~\! dx\leq \irn|x|^{-\beta q}V^q~\!dx~\!.
$$
Since evidently 
$$
\|\Dshalf V\|_2=\|\big(\Dshalf u\big)^*\|_2=\|\Dshalf u\|_2,
$$
we can conclude that
$$
S_q=\frac{\|\Dshalf u\|_2^2}{\|~\!|x|^{-\beta}u\|_q^2}> \frac{\|\Dshalf V\|_2^2}{\|~\!|x|^{-\beta}V\|_q^2}\ge S_q,
$$
and the claim follows via a  contradiction argument.
\medskip

\noindent
{\bf Step 3}. It remains to show that $u$ is strictly decreasing. We recall that $u$ is a minimizer in (\ref{eq:HS_ineq}) and therefore satisfies the integral equation
\begin{equation*}
 u(x)=\irn \frac {|y|^{-\beta q}u^{q-1}(y)}{|x-y|^{n-2s}}\,dy.
\end{equation*}
up to a multiplicative constant.
The right-hand side here is the convolution of two symmetric-decreasing functions, and one of them is strictly symmetric-decreasing. This completes the proof.
\QED

\begin{Remark}
In contrast to the critical case $q=\sstar$, the value of the infimum $S_q$ and its minimizers are not explicitly known even for $s<1$. The regularity, uniqueness and nondegeneracy of the minimizer in this case has been recently proved in \cite{MN20}, see also \cite{BQ, ABQ}. For fractional $s>1$, these properties have never been investigated.
\end{Remark}

\end{document}